\scrollmode % suppress annoying error messages
\documentclass{llncs}
\usepackage{amssymb}
\usepackage{mathrsfs}
\newcommand{\R}{{\mathbb R}}
\newcommand{\A}{{\mathscr A}}

\newcommand{\Z}{{\mathbb Z}}
\newcommand{\Q}{{\mathbb Q}}

\newcommand{\one}{{\mathbf{1}}}

\newcommand{\MMM}{{\mathscr M}}

\renewcommand{\O}{{\mathscr O}}

\newcommand{\fff}{{\mathfrak f}}
\newcommand{\bbb}{{\mathfrak b}}
\newcommand{\aaa}{{\mathfrak a}}
\newcommand{\omitterms}{{\!\!\!\!\phantom{a}^{'}}}
\newcommand{\eee}{{{\mathbf{e}}}}
\newcommand{\eqref}[1]{{(\ref{#1})}}
\newcommand{\sets}[1]{[\![#1]\!]}
\newcommand{\Index}[1]{\|#1\|}
\newcommand{\Sign}{{\mathrm{sgn}}}
\newcommand{\Trace}{{\mathrm{Tr}}}
\newcommand{\rank}{{\mathrm{Rank\ }}}

\begin{document}
\title{Computing special values of partial zeta functions}
\titlerunning{Computing special values}  % abbreviated title (for running head)
%                                     also used for the TOC unless
%                                     \toctitle is used
%
\author{Gautam Chinta\inst{1} \and Paul E. Gunnells\inst{1} \and Robert Sczech\inst{2}}
\authorrunning{G. Chinta, P. E. Gunnells, R. Sczech}   % abbreviated author list (for running head)
%
%%%% modified list of authors for the TOC (add the affiliations)
\tocauthor{Gautam Chinta, Paul E. Gunnells (Columbia University),
Robert Sczech (Rutgers)
}
\institute{Dept. of Mathematics, Columbia University, New York, NY
10027, USA \and Dept. of Mathematics and Computer Science, Rutgers
University, Newark, NJ 07102--1811, USA}

\maketitle              % typeset the title of the contribution

\begin{abstract}
We discuss computation of the special values of partial zeta functions
associated to totally real number fields.  The main tool is the
\emph{Eisenstein cocycle} $\Psi $, a group cocycle for $GL_{n} (\Z )$;
the special values are computed as periods of $\Psi $, and are
expressed in terms of generalized Dedekind sums.  We conclude with
some numerical examples for cubic and quartic fields of small
discriminant.
\end{abstract}
\section{Introduction}
Let $K/\Q $ be a totally real number field of degree $n$ with ring of
integers $\O _{K}$, and let $U\subset \O _{K}^{\times }$ be the
subgroup of totally positive units.  Let $\fff , \bbb \subset \O _{K}$
be relatively prime ideals.  Then the \emph{partial zeta function}
associated to this data is defined by
\[
\zeta _{\fff } (\bbb , s) := \sum _{\aaa \sim \bbb } N (\aaa )^{-s},
\] 
where $\aaa \sim \bbb$ means $\aaa \bbb ^{-1} = (\alpha )$, where
$\alpha $ is a totally positive number in $1+\fff \bbb ^{-1}$.  According to a
classical result of Klingen and Siegel, the special values $\zeta
_{\fff } (\bbb , k)$ are rational for nonpositive integers $k$.
Moreover, the values $\zeta_{\fff } (\bbb , 0)$ are especially
important because of their connection with the Brumer-Stark conjecture
and the Leopoldt conjecture \cite{hayes0,hayes1,hayes2,hayes3,washington}.

In \cite{sczech}, one of us (RS) gave a cohomological interpretation of
these special values by showing that they can be computed in finite
terms as periods of the \emph{Eisenstein cocycle}.
This is a cocycle $\Psi \in H^{n-1} (GL_{n} (\Z ); \MMM )$, where
$\MMM $ is a certain $GL_{n} (\Z )$-module.  Then two of us
(PEG and RS) showed in \cite{dedekind} that the Eisenstein cocycle is an
effectively computable object.  More precisely, using the cocycle one
can express $\zeta _{\fff } (\bbb, k)$ as a finite sum of
\emph{generalized Dedekind sums}, and that the
latter can be effectively computed by a continued-fraction algorithm
that uses a generalization of the classical Dedekind-Rademacher
reciprocity law.

In this note we describe an ongoing project to build a database of
$\zeta_{\fff } (\bbb , 0)$ for various fields $K$ and ideals $\fff,
\bbb $.  We recall the definition of the Eisenstein cocycle and its
relation to the special values (\S\ref{eis.ded}), and discuss the
effective computation of Dedekind sums (\S\ref{diag.uni}).  We
conclude with examples of special values for some fields of degree
3 and 4 (\S\ref{examples}).

\section{Dedekind sums and the Eisenstein cocycle}\label{eis.ded}
\subsection{}\label{gen.dede}
Let $\sigma $ be a square matrix with integral columns
$\sigma _{j}\in \Z ^{n}$ ($j=1,\dots ,n$), and let $L\subset \Z ^{n}$
be a lattice of rank $r\geq 1$.  Let $v\in \Q^{n}$, and let $e\in \Z
^{n}$ with $e_{j}\geq 1$.  Then the \emph{Dedekind sum} $S$ associated to
the data $(L,\sigma ,e,v)$ is defined by
\begin{equation}\label{Dede.sum}
S = S (L,\sigma ,e,v) := \sum _{x\in L} \omitterms \eee (\langle
x,v\rangle )\frac{\det \sigma}{\langle x,\sigma _{1}\rangle
^{e_{1}}\cdots \langle x,\sigma _{n}\rangle ^{e_{n}}}.
\end{equation}
Here $\langle x,y\rangle := \sum x_{i}y_{i}$ is the usual scalar
product on $\R ^{n}$, $\eee (t)$ is the character $\exp (2\pi it)$,
and the prime next to the summation means to omit terms for which the
denominator vanishes.  The series \eqref{Dede.sum} converges
absolutely if all $e_{j}>1$, but may only converge conditionally if $e_{j}=1$ for
some $j$. In this latter case we can define the sum by the
\emph{$Q$-limit}
\begin{equation}\label{qlim}
\sum _{x\in L}\omitterms a (x)\Bigr|_{Q} := \lim _{t\rightarrow \infty
}\left(\sum_{|Q (x)|<t}\!\!\!\!\omitterms  a (x)  \right),
\end{equation}
where $Q$ is any finite product of real-valued linear forms on $\R
^{n}$ that doesn't vanish on $\Q ^{n}\smallsetminus \{0 \}$.  One can
precisely determine how the value of \eqref{Dede.sum} depends on $Q$
(\cite[Thm. 7]{sczech}).  The sum $S$ is always a rational number
times a power of $2\pi i$.

\subsection{}\label{eis.coc}
We recall now the definition of the Eisenstein cocycle $\Psi $ and its
relationship with the special values $\zeta _{\fff } (\bbb , k)$.  For
simplicity, we describe only material necessary to compute the special
value at $k=0$, and refer to \cite{sczech,dedekind} for other $k$.

Let $\A = (A_{1},\dots ,A_{n}) \in (GL_{n} (\R ))^{n}$ be an $n$-tuple of matrices.  For an $n$-tuple $d = (d_{1},\dots ,d_{n})$ of integers
$1\leq d_{i}\leq n$, let $\A (d)\subseteq \R ^{n}$ be the subspace
generated by all columns $A_{ij}$ such that $j<d_{i}$.  (Here $A_{ij}$
denotes the $j$th column of the matrix $A_{i}$.)  Writing $\A (d)^{\perp}$ for
the orthogonal complement of $\A (d)$ in $\R^{n}$, we let 
\begin{equation}\label{one}
X (d) = \A (d)^{\perp } \smallsetminus \bigcup_{i=1}^{n} \sigma
_{i}^{\perp }, \quad \hbox{where $\sigma _{i} = A_{id_{i}}$.}
\end{equation}
The $n$-tuple $\A $ determines a decomposition of $\R ^{n}
\smallsetminus \{0 \}$ into linear strata
\begin{equation}
\bigsqcup_{d\in D} X (d),
\end{equation}
indexed by the finite set 
\[
D = D (\A ) = \{d \mid X (d) \not = \varnothing  \}.
\]
Associated to this decomposition is a collection of rational functions
$\psi (\A )$ on $\R ^{n}\smallsetminus \{0 \}$, defined by
\[
\psi (\A ) (x) = \frac{\det (\sigma _{1},\dots ,\sigma _{n})}{\langle
x,\sigma _{1}\rangle \cdots \langle x,\sigma _{n}\rangle },\quad
\hbox{if $x\in X (d)$.}
\]
Note that $\psi (\A ) (x)$ is well-defined by the construction of $X
(d)$.

Let $v\in \R ^{n}$, and let $Q$ be defined as in \S\ref{gen.dede}.
Then the \emph{Eisenstein cocycle} $\Psi $ is defined as
\[
\Psi = \Psi (\A ) (Q,v) := (2\pi i)^{-n}\sum _{x\in \Z ^{n}} \eee (\langle x,v\rangle )
\psi (\A ) (x)\Bigr|_{Q}.
\]
One can show that $\Psi $ is a homogeneous $(n-1)$-cocycle
for $GL_{n} (\Z )$.  Furthermore, we can express $\Psi $ in terms of
Dedekind sums 
\begin{equation}\label{dedsums}
\Psi (\A ) (Q,v) = (2\pi i)^{-n}\sum _{d\in D} S (L(d),\sigma ,\one
,v)\Bigr|_{Q},
\end{equation}
where $\sigma $ is the matrix with columns $A_{id_{i}},$  ($i=1,\dots ,n$),
$L(d)$ is the lattice $\A (d)^\perp \cap \Z ^{n}$, and $\one $ is the
vector $(1,\dots ,1)$.

\subsection{}
Now we describe how $\Psi $ can be used to compute special values.
Let $W$ be a $\Z $-basis for the fractional ideal $\fff \bbb ^{-1} =
\sum \Z W_{j}$, and let $W^{*}$ be the dual basis with respect to the
trace form.  Via the $n$ real embeddings $\tau _{i}$, $i=1,\dots ,n$, any $x\in K$ determines a row vector $(\tau_{1} (x),\dots
,\tau_{n} (x))$.  Hence we may identify $W$ with a matrix in
$GL_{n} (\R )$: the $j$th row of this matrix is the image of the $j$th
basis element of $W$.  Let
\[
Q (X)= \prod _{i}\sum _{j}X_{j}(\tau _{i} ( W_{j}^{*})),
\]
and let $v\in \Q ^{n}$ be defined by $v_{j} =
\Trace (W^{*}_{j})$.

Let $\nu = n-1$, and let $\varepsilon _{1},\dots ,\varepsilon _{\nu }$
be a basis for the totally positive units $U$.  Using the regular
representation $\rho $ with respect to the basis $W$, we identify the units
$\varepsilon _{j}$ with elements $A_{j} = \rho (\varepsilon_{j} )^{t}
\in GL_{n} (\Z )$.  Using the bar notation
\[
[A_{1}|\cdots | A_\nu ] := (1,A_{1},A_{1}A_{2},\dots ,A_{1}\cdots
A_{\nu })\in (GL_{n} (\Z ))^{n},
\]
we have the following proposition expressing the zeta values in terms
of the Eisenstein cocycle:
\begin{proposition}\label{zeta-values}
\cite{sczech,dedekind}
Let $U_{\fff }$ be the subgroup $U\cap (1+\fff )$, and let $\pi $ run
through all permutations of \hbox{$\{1,\dots ,\nu \}$}.  Then 
\[
\zeta_{\fff } (\bbb ,0) = \eta \sum _{\varepsilon \in U/U_{\fff }}
\sum _{\pi } \Sign (\pi) \Psi ([A_{\pi (1) }|\cdots | A_{\pi (\nu )}])
(Q,\rho (\varepsilon )^{t}v).
\]
Here $\eta = \pm
1$ is defined by 
\[
\eta = (-1)^{\nu }\Sign (\det W) \Sign (R), 
\]
where $R = \det (\log \tau_{j}(\varepsilon _{i}))$, $1\leq i,j\leq \nu $.
\end{proposition}

\section{Diagonality and unimodularity}\label{diag.uni}

\subsection{}
We define the \emph{rank} of $S=S (L,\sigma ,e,v)$ to be the rank of
the lattice $L$.  It is easy to see that after a $GL_n(\Q)$
transformation, we may assume that $L$ is the sublattice $Z ^{\ell }$
spanned by the first $\ell $ standard basis vectors, where $\ell $ is
the rank of $L.$ Furthermore, by multiplying by an appropriate
rational factor, permuting columns and repeating columns if necessary,
we may assume the pair $(Z ^{\ell },\sigma)$ satisfies the following
conditions:
\begin{enumerate}
\item [(i)] For each column $\sigma_{j} $, the vector of the first
$\ell $ components of $\sigma _{j}$ is primitive and integral.
\item [(ii)] If two columns of $\sigma$ induce proportional linear forms on
$Z ^{\ell }$, then these two linear forms coincide on $Z ^{\ell }$,
and are adjacent columns of $\sigma$.  
\item [(iii)]  The vector $e=\one $.
\end{enumerate}

Let $S (Z ^{\ell },\sigma ,\one ,v)$ be a Dedekind sum satisfying the
three conditions above.  Let $\pi \colon \R ^{N} \rightarrow \R ^{\ell
}$ be the projection on the first $\ell $ components, and let $\pi
(\sigma )$ be the $\ell \times n $ matrix with columns $\pi (\sigma
_{i})$.

\begin{definition}\label{index}
Let $M (\sigma )$ be the set of maximal minors of $\pi (\sigma )$.
Then the \emph{index} of $S$, denoted $\Index{S}$, is defined to be
\[
\max_{\tau \in M (\sigma )} |\det \tau |.
\]
A Dedekind sum is \emph{unimodular} if $\Index{S} = 1$.
 
\end{definition}

\subsection{}
Now define a partition 
\begin{equation}\label{part}
\sets{n} = \bigsqcup_{k=1}^{s} I_{k}, \quad
\hbox{$\ell \leq s\leq n$}
\end{equation}
as follows.  Put
\[
i,j\in I_{k} \quad \hbox{if and only if}\quad \pi
(\sigma _{i}) = \pi (\sigma _{j}).
\]
In other words, two elements of $\sets{n}$ are in the same set of the
partition if the corresponding columns of $\sigma $ induce
the same linear form on $Z ^{\ell }$.

Let $p_{k}=\#I_{k}$.   

\begin{definition}\label{diag}
The vector $p(S) = (p_{1},\dots ,p_{s})$ is called the \emph{type} of
$S$. A Dedekind sum is called \emph{diagonal} if $p (S)$ has length
$\ell$.
\end{definition}

\subsection{}
The virtue of diagonality is that a diagonal Dedekind sum $S$ may be
evaluated as a finite sum of products of generalized Bernoulli
polynomials.  Furthermore, the number of terms in this finite sum is
the index of $S$.   Hence diagonal and unimodular Dedekind sums can be 
evaluated very rapidly.

In general, the Dedekind sums in \eqref{dedsums} aren't diagonal.
However, we have the following theorem, which is the main result of
\cite{dedekind}:

\begin{theorem}\label{mnthm}
\cite{dedekind}
Every Dedekind sum $S (L,\sigma ,e,v)$ can be expressed as a finite
rational linear combination of unimodular diagonal sums.  If 
$n$, $\rank L$, and $e$ are fixed, then this expression can be
computed in time polynomial in $\log\Index{S}$.  Moreover, the number of
terms in this expression is bounded by a polynomial in $\log\Index{S}$.
\end{theorem}

The key ingredient in the proof of Theorem \ref{mnthm} is a
``reciprocity law'' for higher-dimensional Dedekind sums.
For any nonzero point $v\in \R ^{n}$, let $v^{\perp }$ be the
hyperplane $\{x \mid \langle v,x\rangle =0 \}$.  Let $Q$ be a finite
product of real-valued linear forms on $\R ^{n}$ that do not vanish on
$\Q ^{n}\smallsetminus \{0 \}$.

\begin{proposition}\label{recip}
Let $\sigma _{0},\dots ,\sigma _{n}\in \Z ^{n}$ be nonzero.  For
$j=0,\dots ,n$, let $\sigma ^{j}$ be the matrix with columns $\sigma
_{0},\dots ,\hat \sigma _{j},\dots ,\sigma _{n}$.  Fix a lattice
$L\subseteq \Z ^{n}$, and assume $e = \one$.  Then for any $v\in \R
^{n}$, we have the following identity among Dedekind sums:
\begin{equation}\label{recip.law}
\sum _{j=0}^{n} (-1)^{j} S (L,\sigma ^{j},\one,v )\Bigr|_{Q} =
\sum _{j=0}^{n} (-1)^{j} S (L\cap \sigma _{j}^{\perp },\sigma
^{j},\one,v )\Bigr|_{Q}. 
\end{equation}

\end{proposition}

We refer to \cite{dedekind} for proofs of the above statements.  Here,
in the following two sections, we show how Theorem \ref{mnthm} is
applied with a rank $2$ example.   For simplicity we ignore issues of
convergence, and merely remark that all of our manipulations with sums
are compatible with the $Q$-limit process \eqref{qlim}.

\subsection{}\label{previous}
Let $L$ be the lattice $Z^2$.  Let
\[
\sigma=\left( \begin{array}{ccc}
      1 & 0 & 1 \\
      0 & 1 & 2 \\
      0 & 0& 1\end{array} \right),\quad e=(1,1,2), \quad \hbox{and $v=(0,0,0)$.}
\]
Hence $S(L,\sigma,e,v)$ denotes the absolutely convergent sum
\[
\sum_{(x,y)\in\Z^2} \!\!\!\!\omitterms\ \ \frac{1}{xy(x+2y)^2},
\]
where the prime on the summation indicates that we omit the terms
$(x,y)$ for which $x$, $y$ or $x+2y$ vanish.

This sum isn't diagonal, since $\sigma $ induces $3$ different
linear forms on $L$ instead of $2$.  To diagonalize $S$, 
we begin with the identity of rational functions
\begin{equation}
\label{pf}
\frac{1}{xy(x+2y)^2}=\frac{1}{y(x+2y)^3}+\frac{2}{x(x+2y)^3}.
\end{equation}
This is true provided none of the denominators vanishes.  The
numerators of the functions on the right come from expressing the
third column of $\sigma$ as a linear combination of the first two:
$$(1,2)^{t}=1\cdot(1,0)^{t}+2\cdot(0,1)^{t}.$$

We want to sum both sides of \eqref{pf} over pairs $(x,y)\in \Z ^{2}$
to obtain an identity among Dedekind sums of the form
\begin{equation}
\label{formalId}
\sum_{(x,y)\in\Z^2} \!\!\!\!\omitterms\ \ \frac{1}{xy(x+2y)^2}=
\sum_{(x,y)\in\Z^2} \!\!\!\!\omitterms\ \ \frac{1}{y(x+2y)^3}+
\sum_{(x,y)\in\Z^2} \!\!\!\!\omitterms\ \ \frac{2}{x(x+2y)^3}.
\end{equation}
However, as written \eqref{formalId} is incorrect.  The identity
\eqref{pf} only holds if none of $x$, $y$, or $x+2y$ vanish, but the
sums on the right of \eqref{formalId} include some of these terms (for
instance, the first sum on the right of \eqref{formalId}
contains terms $(x,y)$ with $x=0$).  We account for this by
subtracting two rank $1$ Dedekind sums from the right of
\eqref{formalId} as ``correction terms'': 
\begin{eqnarray}\nonumber
\sum_{(x,y)\in\Z^2} \!\!\!\!\omitterms\ \ \frac{1}{xy(x+2y)^2}&=&
\sum_{(x,y)\in\Z^2} \!\!\!\!\omitterms\ \ \frac{1}{y(x+2y)^3}
  + \sum_{(x,y)\in\Z^2} \!\!\!\!\omitterms\ \ \frac{2}{x(x+2y)^3}
 \\ \nonumber
&\qquad\qquad & - \mathop{\sum_{(x,y)\in\Z^2}}_{x=0}\!\!\!\!\omitterms \ \ \frac{1}{y(x+2y)^3} -
\mathop{\sum_{(x,y)\in\Z^2}}_{y=0}\!\!\!\!\omitterms \ \ \frac{1}{x(x+2y)^3}\\ \nonumber
\label{recip.law.ex}\nonumber
&=&
\sum_{(x,y)\in\Z^2} \!\!\!\!\omitterms\ \ \frac{1}{y(x+2y)^3} + \sum_{(x,y)\in\Z^2} \!\!\!\!\omitterms\ \ \frac{2}{x(x+2y)^3} \\ 
&\qquad\qquad & -
\sum_{y\in\Z}\omitterms\frac{1}{8y^4} \qquad\qquad -
\sum_{x\in \Z}\omitterms\frac{2}{x^4}.
\end{eqnarray}
This equation is precisely an instance of the reciprocity law
(Proposition \ref{recip}).  The three rank $2$ sums are the left of
\eqref{recip.law}, and the two rank $1$ sums are the right of
\eqref{recip.law} (one rank $1$ sum in \eqref{recip.law} vanishes
identically).  Note that all of the sums on the right of
\eqref{recip.law.ex} are now diagonal.  

To diagonalize a general Dedekind sum $S (L,\sigma ,\one ,v)$, one
considers the configuration $C\subset \R ^{n}$ of linear subspaces
consisting of $(L\otimes \R )^{\perp}$ and the spaces generated by the
points $\sigma _{1},\dots ,\sigma _{n}$.  One shows by investigating
the geometry of $C$ that a point $\sigma _{0}$ can be found such that
when Proposition \ref{recip} is applied with the tuple $(\sigma
_{0},\dots ,\sigma _{n})$, the resulting Dedekind sums are ``closer''
to diagonality in a certain sense.  It may take several applications
of Proposition \ref{recip} to express a Dedekind sum as a linear
combination of diagonal sums.

\subsection{}\label{where.rho.is}
The second rank two sum on the right of \eqref{recip.law.ex}
has index two.  We will show how to make this sum unimodular.  Write
$$S(\Z^2,\tau,(1,3),v')=\sum_{(x,y)\in\Z^2}\!\!\!\! \omitterms\ \ \frac{1}{x(x+2y)^3}$$
where
$$\tau=(\tau_1,\tau_2)=\left( \begin{array}{cc}
      1 & 1 \\
      0 & 2 \end{array} \right)\quad \hbox{and}\quad v'=(0,0).$$

Let $\rho$ be the column vector $(0,1)^{t}$.  We apply Proposition \ref{recip}
to the triple $(\rho,\tau_1,\tau_2)$:
\begin{eqnarray}\label{unimod.ex}\nonumber
\sum_{(x,y)\in\Z^2} \!\!\!\! \omitterms\ \ \frac{2}{x(x+2y)^3}&=&
\sum_{(x,y)\in\Z^2} \!\!\!\! \omitterms\ \ \frac{-1}{y(x+2y)^3}+\sum_{(x,y)\in\Z^2} \!\!\!\! \omitterms\ \ \frac{1}{x(x+2y)^2y}\\ 
&\qquad\qquad & +\sum_{y\in\Z}\omitterms \frac{1}{8y^4} \qquad\qquad +
\sum_{x\in \Z}\omitterms \frac{2}{x^4}.
\end{eqnarray}
Now all the terms on the right of \eqref{unimod.ex} are diagonal and
unimodular except for the second rank two sum.  In fact, this sum is no
longer diagonal.  However, one further application of Proposition
\ref{recip} as in \S\ref{previous} will make the third sum diagonal
and unimodular.  Hence we will have succeeded in expressing the
original sum as a finite linear combination of
diagonal, unimodular Dedekind sums.

In general, one must be able to construct the vector $\rho $ as
above.  An easy argument using Minkowski's Theorem from the geometry
of numbers guarantees the existence of $\rho $ \cite{ash.rudolph}.
To construct $\rho $ in practice, one may use $LLL$-reduction of the
lattice spanned by the \emph{rows} of $\sigma $ and \cite[Conjecture
3.9]{experimental}.

\section{Examples}\label{examples}
Here we present some numerical examples.  For simplicity we compute
$\zeta = \zeta _{\fff } (\bbb ,0)$, where $\fff = N\O_{K}$ for various
rational integers $N$, and $\bbb = \O_K$.  These fields are the
first entries in the tables of totally real fields with small
discriminant, available from \cite{bordeaux}.

\bigskip
\noindent \textsl{Cubic fields}

\begin{itemize}
\item [$\bullet$] $K=\Q (\theta )$, where $\theta ^3+\theta ^2-2\theta
-1=0$ (discriminant $49$).

\begin{center}
\begin{tabular}{|c|r||c|r||c|r||c|r||c|r||c|r||c|r||c|r|}
\hline
$N$&$N\cdot \zeta$&$N$&$N\cdot \zeta$&$N$&$N\cdot \zeta$&$N$&$N\cdot \zeta$&$
N$&$N\cdot \zeta$&$N$&$N\cdot \zeta$&$N$&$N\cdot \zeta$&$N$&$N\cdot \zeta$\\
\hline\hline
&&$4$&$1$&$7$&$2$&$10$&$16$&$13$&$-1$&$16$&$-9$&$19$&$18$&$22$&$-8$\\
$2$&$0$&$5$&$2$&$8$&$3$&$11$&$-10$&$14$&$-24$&$17$&$-26$&$20$&$19$&$23$&$10$\\
$3$&$2$&$6$&$8$&$9$&$2$&$12$&$5$&$15$&$14$&$18$&$8$&$21$&$6$&$24$&$23$\\
\hline
\end{tabular}
\end{center}

\newpage

\item [$\bullet$] $K=\Q (\theta )$, where $\theta ^3-3\theta
-1=0$ (discriminant $81$).

\begin{center}
\begin{tabular}{|c|r||c|r||c|r||c|r||c|r||c|r||c|r||c|r|}
\hline
$N$&$N\cdot \zeta$&$N$&$N\cdot \zeta$&$N$&$N\cdot \zeta$&$N$&$N\cdot \zeta$&$
N$&$N\cdot \zeta$&$N$&$N\cdot \zeta$&$N$&$N\cdot \zeta$&$N$&$N\cdot \zeta$\\
\hline\hline
&&$4$&$1$&$7$&$14$&$10$&$-16$&$13$&$2$&$16$&$-13$&$19$&$11$&$22$&$72$\\
$2$&$0$&$5$&$-2$&$8$&$7$&$11$&$-6$&$14$&$0$&$17$&$13$&$20$&$11$&$23$&$62$\\
$3$&$2/3$&$6$&$8/3$&$9$&$2/3$&$12$&$11/3$&$15$&$14/3$&$18$&$-64/3$&$21$&$14/3$&$24$&$29/3$\\
\hline
\end{tabular}
\end{center}

\item [$\bullet$] $K=\Q (\theta )$, where $\theta ^3+\theta ^2-3\theta
-1= 0$ (discriminant $148$).

\begin{center}
\begin{tabular}{|c|r||c|r||c|r||c|r||c|r||c|r||c|r||c|r|}
\hline
$N$&$N\cdot \zeta$&$N$&$N\cdot \zeta$&$N$&$N\cdot \zeta$&$N$&$N\cdot \zeta$&$
N$&$N\cdot \zeta$&$N$&$N\cdot \zeta$&$N$&$N\cdot \zeta$&$N$&$N\cdot \zeta$\\
\hline\hline
&&$4$&$1$&$7$&$2$&$10$&$-2$&$13$&$-22$&$16$&$7$&$19$&$82$&$22$&$68$\\
$2$&$0$&$5$&$-4$&$8$&$3$&$11$&$-18$&$14$&$-20$&$17$&$100$&$20$&$4$&$23$&$12$\\
$3$&$2$&$6$&$4$&$9$&$-10$&$12$&$5$&$15$&$42$&$18$&$-32$&$21$&$-78$&$24$&$23$\\
\hline
\end{tabular}
\end{center}

\item [$\bullet$] $K=\Q (\theta )$, where $\theta ^3-\theta ^2-4\theta -1= 0$
(discriminant $169$).

\begin{center}
\begin{tabular}{|c|r||c|r||c|r||c|r||c|r||c|r||c|r||c|r|}
\hline
$N$&$N\cdot \zeta$&$N$&$N\cdot \zeta$&$N$&$N\cdot \zeta$&$N$&$N\cdot \zeta$&$
N$&$N\cdot \zeta$&$N$&$N\cdot \zeta$&$N$&$N\cdot \zeta$&$N$&$N\cdot \zeta$\\
\hline\hline
&&$4$&$3$&$7$&$6$&$10$&$8$&$13$&$2$&$16$&$25$&$19$&$-238$&$22$&$160$\\
$2$&$0$&$5$&$1$&$8$&$5$&$11$&$-6$&$14$&$56$&$17$&$6$&$20$&$-8$&$23$&$386$\\
$3$&$2$&$6$&$8$&$9$&$26$&$12$&$11$&$15$&$17$&$18$&$52$&$21$&$-10$&$24$&$89$\\
\hline
\end{tabular}
\end{center}

\item [$\bullet$] $K=\Q (\theta )$, where $\theta ^3-4\theta -1= 0$
(discriminant $229$).

\begin{center}
\begin{tabular}{|c|r||c|r||c|r||c|r||c|r||c|r||c|r||c|r|}
\hline
$N$&$N\cdot \zeta$&$N$&$N\cdot \zeta$&$N$&$N\cdot \zeta$&$N$&$N\cdot \zeta$&$
N$&$N\cdot \zeta$&$N$&$N\cdot \zeta$&$N$&$N\cdot \zeta$&$N$&$N\cdot \zeta$\\
\hline\hline
&&$4$&$2$&$7$&$18$&$10$&$-24$&$13$&$60$&$16$&$19$&$19$&$26$&$22$&$8$\\
$2$&$0$&$5$&$-2$&$8$&$7$&$11$&$-4$&$14$&$24$&$17$&$8$&$20$&$38$&$23$&$202$\\
$3$&$2$&$6$&$0$&$9$&$-10$&$12$&$18$&$15$&$56$&$18$&$-16$&$21$&$96$&$24$&$51$\\
\hline
\end{tabular}
\end{center}

\item [$\bullet$] $K=\Q (\theta )$, where $\theta ^3-\theta ^2-4\theta +3= 0$
(discriminant $257$).

\begin{center}
\begin{tabular}{|c|r||c|r||c|r||c|r||c|r||c|r||c|r||c|r|}
\hline
$N$&$N\cdot \zeta$&$N$&$N\cdot \zeta$&$N$&$N\cdot \zeta$&$N$&$N\cdot \zeta$&$
N$&$N\cdot \zeta$&$N$&$N\cdot \zeta$&$N$&$N\cdot \zeta$&$N$&$N\cdot \zeta$\\
\hline\hline
&&$4$&$6$&$7$&$16$&$10$&$4$&$13$&$10$&$16$&$5$&$19$&$42$&$22$&$72$\\
$2$&$0$&$5$&$3$&$8$&$9$&$11$&$-6$&$14$&$0$&$17$&$-26$&$20$&$12$&$23$&$-112$\\
$3$&$2$&$6$&$8$&$9$&$-4$&$12$&$8$&$15$&$21$&$18$&$-16$&$21$&$-12$&$24$&$58$\\
\hline
\end{tabular}
\end{center}

\end{itemize}

\bigskip
\noindent \textsl{Quartic fields}

\begin{itemize}
\item [$\bullet$] $K=\Q (\theta )$, where $\theta ^4 - \theta ^3 -
3\theta ^2 + \theta + 1 = 0$ (discriminant $725$).
\begin{center}
\begin{tabular}{|c|r||c|r||c|r||c|r||c|r||c|r||c|r||c|r|}
\hline
$N$&$N\cdot \zeta$&$N$&$N\cdot \zeta$&$N$&$N\cdot \zeta$&$N$&$N\cdot \zeta$&$
N$&$N\cdot \zeta$&$N$&$N\cdot \zeta$&$N$&$N\cdot \zeta$&$N$&$N\cdot \zeta$\\
\hline\hline
&&$4$&$1$&$7$&$-4$&$10$&$0$&$13$&$-20$&$16$&$35$&$19$&$-32$&$22$&$32$\\
$2$&$0$&$5$&$4$&$8$&$7$&$11$&$-2$&$14$&$-16$&$17$&$-92$&$20$&$0$&$23$&$-12$\\
$3$&$4$&$6$&$0$&$9$&$-44$&$12$&$0$&$15$&$0$&$18$&$-320$&$21$&$-30$&$24$&$-84$\\
\hline
\end{tabular}
\end{center}

\newpage

\item [$\bullet$] $K=\Q (\theta )$, where $\theta ^4-\theta ^3-4\theta ^2+4\theta +1= 0$ (discriminant $1125$). 
\begin{center}
\begin{tabular}{|c|r||c|r||c|r||c|r||c|r||c|r||c|r||c|r|}
\hline
$N$&$N\cdot \zeta$&$N$&$N\cdot \zeta$&$N$&$N\cdot \zeta$&$N$&$N\cdot \zeta$&$
N$&$N\cdot \zeta$&$N$&$N\cdot \zeta$&$N$&$N\cdot \zeta$&$N$&$N\cdot \zeta$\\
\hline\hline
&&$3$&$4/5$&$5$&$4/3$&$7$&$124/3$&$9$&$-116/5$&$11$&$72$&$13$&$676/3$&$15$&$4$\\
$2$&$0$&$4$&$2$&$6$&$0$&$8$&$7$&$10$&$128$&$12$&$0$&$14$&$144$&$16$&$51$\\

\hline
\end{tabular}
\end{center}

\item [$\bullet$] $K=\Q (\theta )$, where $\theta ^4-6\theta ^2+4=0$ (discriminant $1600$).
\begin{center}
\begin{tabular}{|c|r||c|r||c|r||c|r||c|r||c|r||c|r||c|r|}
\hline
$N$&$N\cdot \zeta$&$N$&$N\cdot \zeta$&$N$&$N\cdot \zeta$&$N$&$N\cdot \zeta$&$
N$&$N\cdot \zeta$&$N$&$N\cdot \zeta$&$N$&$N\cdot \zeta$&$N$&$N\cdot \zeta$\\
\hline\hline
&&$4$&$1/2$&$7$&$4$&$10$&$4$&$13$&$104$&$16$&$45/2$&$19$&$224$&$22$&$136$\\
$2$&$1$&$5$&$4$&$8$&$5/2$&$11$&$-4$&$14$&$-16$&$17$&$-84$&$20$&$-9$&$23$&$60$\\
$3$&$4$&$6$&$0$&$9$&$16$&$12$&$24$&$15$&$6$&$18$&$128$&$21$&$222$&$24$&$15$\\
\hline
\end{tabular}
\end{center}

\item [$\bullet$] $K=\Q (\theta )$, where $\theta ^4-4\theta ^2-\theta +1=0$ (discriminant $1957$).
\begin{center}
\begin{tabular}{|c|r||c|r||c|r||c|r||c|r||c|r||c|r|}
\hline
$N$&$N\cdot \zeta$&$N$&$N\cdot \zeta$&$N$&$N\cdot \zeta$&$N$&$N\cdot \zeta$&$
N$&$N\cdot \zeta$&$N$&$N\cdot \zeta$&$N$&$N\cdot \zeta$\\
\hline\hline
&&$3$&$4$&$5$&$20$&$7$&$-8$&$9$&$52$&$11$&$4$&$13$&$500$\\
$2$&$0$&$4$&$1$&$6$&$0$&$8$&$11$&$10$&$0$&$12$&$3$&$14$&$-8$\\
\hline
\end{tabular}
\end{center}

\item [$\bullet$] $K=\Q (\theta )$, where $\theta ^4-5\theta ^2+5=0$ (discriminant $2000$).
\begin{center}
\begin{tabular}{|c|r||c|r||c|r||c|r||c|r||c|r||c|r||c|r|}
\hline
$N$&$N\cdot \zeta$&$N$&$N\cdot \zeta$&$N$&$N\cdot \zeta$&$N$&$N\cdot \zeta$&$
N$&$N\cdot \zeta$&$N$&$N\cdot \zeta$&$N$&$N\cdot \zeta$&$N$&$N\cdot \zeta$\\
\hline\hline
&&$4$&$11/5$&$7$&$52$&$10$&$0$&$13$&$-32$&$16$&$31/10$&$19$&$412/5$&$22$&$568$\\
$2$&$2/5$&$5$&$0$&$8$&$71/10$&$11$&$4$&$14$&$16$&$17$&$296$&$20$&$30$&$23$&$-148$\\
$3$&$4$&$6$&$0$&$9$&$4$&$12$&$12$&$15$&$340$&$18$&$-112$&$21$&$428$&$24$&$117$\\
\hline
\end{tabular}
\end{center}

\item [$\bullet$] $K=\Q (\theta )$, where $\theta ^4-4\theta ^2+2=0$ (discriminant $2048$).
\begin{center}
\begin{tabular}{|c|r||c|r||c|r||c|r||c|r||c|r||c|r|}
\hline
$N$&$N\cdot \zeta$&$N$&$N\cdot \zeta$&$N$&$N\cdot \zeta$&$N$&$N\cdot \zeta$&$
N$&$N\cdot \zeta$&$N$&$N\cdot \zeta$&$N$&$N\cdot \zeta$\\
\hline\hline
&&$3$&$4$&$5$&$-28$&$7$&$8$&$9$&$-20$&$11$&$68$&$13$&$-52$\\
$2$&$1/2$&$4$&$1/4$&$6$&$-10$&$8$&$9/4$&$10$&$19$&$12$&$-41/2$&$14$&$36$\\
\hline
\end{tabular}
\end{center}
\end{itemize}

\bibliographystyle{amsplain}
\bibliography{deadants}

\end{document}